\documentclass{article}
\usepackage{graphics}

\newtheorem{theorem}{Theorem}[section]

\newtheorem{conj}[theorem]{Conjecture}
\newtheorem{prob}[theorem]{Problem}
\newtheorem{assumption}[theorem]{Assumption}

\title{Genera of Vertex Operator Algebras and \\
three-dimensional Topological Quantum Field Theories}

\author{Gerald H\"ohn
\thanks{Mathematisches Institut, Eckerstra\ss{}e~1, 79104 Freiburg, Germany, \hfill\newline
{\tt gerald@mathematik.uni-freiburg.de} }}

\date{}

\def\H{{\bf H}}
\def\C{{\bf C}}
\def\R{{\bf R}}

\def\Z{{\bf Z}}

\def\Q{{\bf Q}}

\def\qed{\hfill\framebox[3mm][t1]{\phantom{x}}}

\begin{document}

\maketitle


\begin{abstract}
The notion of the {\it genus\/} of a quadratic form is generalized to vertex
operator algebras. We define it as the modular braided tensor category
associated to a suitable vertex operator algebra together with the central 
charge. Statements similar as known for quadratic forms are formulated.

We further explain how extension problems for vertex operator algebras can
be described in terms of the associated modular braided tensor category.
\end{abstract}

\maketitle


\section{Introduction}

This note is a write-up of talks given at the workshop
``Vertex Operator Algebras in Mathematics and Physics''
at the Fields Institute in Toronto in October 2000 and at 
some other occasions before. Its purpose is to explain how results 
in low dimensional
quantum field theory obtained from different perspectives and 
with different motivations can be put into a uniform picture by 
generalizing the notions of discriminant forms and genera from integral
quadratic forms to vertex operator algebras, indicating a rich
underlying arithmetical theory.

\medskip 

The relation between two-dimensional conformal field theory and 
three-dimen\-sional topological quantum field theory has been quite known 
{}from the beginning. For example, Witten relates in his well-known
paper~\cite{Wi-tqft} the  Wess-Zumino-Witten models with 
Chern-Simons theories and the Jones polynomial of knots.
The basic idea is to use the monodromy properties of correlation functions 
on surfaces to obtain invariants for knots and three-dimensional manifolds.

Being not overly precise, one can say
that a quantum field theory is a functor from some bordism category
of oriented manifolds to some category of linear spaces satisfying some basic gluing
axioms.

\newpage
\vspace{12mm}
\includegraphics{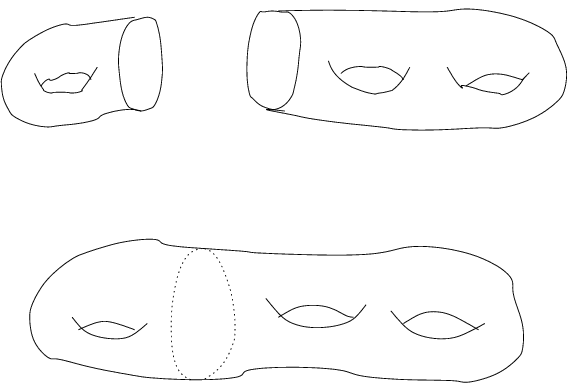}

\vspace{-45mm}
$\ \ \ \  X   \qquad  \qquad  \qquad  \qquad Y $

\vspace{3mm}
$ \ \ \ \ \ \ \ \ \ \ \ \ \ \ \ \ \ \ \ \ \ \ \ \ \ \ \ \ \ \ \ \ \ \
 \ \ \ \ \ \ \ \ \ \ \ \ \ \ \ \ \ \ \ \ \
\stackrel{\mathcal{F}}{\longrightarrow} \ \
\mathcal{F}(X)\otimes  \mathcal{F}(Y) \in  \mathcal{F}(W)\otimes
  \mathcal{F}(W)^*$

\vspace{8mm}
\ \ \ \ \ \ \ \ \  
\hbox{\footnotesize $\ \partial X=W=-\partial Y$}\ \ \ \ 
$ \Downarrow \hbox{gluing}\ \ \ \ \ \ \ \ \ \ \  \ \ \ \ \  \ \ \ \ \  \ 
 \ \ \ \ \ \ \ 
\downarrow  
  \qquad  \qquad  \downarrow  \langle \,.\,,\,.\,\rangle\phantom{|}_{\rm nat}$

\vspace{9mm}
$ \ \ \ \ \ \ \ \ \ \ \ \ \ \ \ \ \ \ \ \ \ \ \ \ \ \ \ \ \ \ \ \ \ \
 \ \ \ \ \ \ \ \ \ \ \ \ \ \ \ \ \ 
\stackrel{\mathcal{F}}{\longrightarrow}
 \ \ \ \ \ \ \ \ \ \ \ \ \ 
\mathcal{F}(Z) \in  \mathcal{F}(\emptyset)\cong\C$

\vspace{8mm}

$\qquad\qquad  Z=X\coprod\nolimits_W Y$

\vspace{5mm}

For two-dimensional conformal field theories, 
one takes two-dimensional surfaces with a conformal
structure on it and one is forced to take infinite dimensional linear spaces.
Using instead Riemannian surfaces (which have a naturally defined 
conformal structure), one reaches the more restricted but essential 
notion of chiral algebras.
For three-dimensional topological quantum field theories
(3d-TQFTs), one takes three-dimensional topological manifolds 
(together with extra structures like framed links) and
only finite dimensional linear spaces.

It took some time before precise and efficient mathematical formulations for
both structures, originating partly from physics, emerged.
Chiral algebras have been axiomatized by
the algebraic structure of {\it vertex algebras\/}; 
see~\cite{Bo-ur,FLM,Se-cft,Hu-Geo2,Kac-VOA,Frenkel,BD} for some approaches. 
3d-TQFTs can be constructed from {\it modular braided tensor categories\/};
see~\cite{MoSei,At-knot,Turaev,BaKi}.

\smallskip

In this paper, we will restrict ourselves mainly to vertex operator algebras (VOAs);
these are vertex algebras which have a Virasoro element and 
we assume that they are non-negatively graded;
we consider only simple and rational VOAs and we make some further strong assumptions. 
There are many different kinds of examples of VOAs known, 
but it is quite unclear yet how a general structure theory for them may look like.

An important class of VOAs can be constructed from positive definite even integral 
lattices. For integral lattices, there is a deeply developed structure theory available 
leading to invariants like rational equivalence, the genus, or the spinor genus.
The first purpose of this paper is to explain how the notion of the {\it genus\/} 
can be generalized in a natural way from lattices to VOAs by using modular braided 
tensor categories so that the definition is compatible with the notion of the genus 
for lattices. Secondly, we show that new VOAs can be constructed from given
ones only by using the associated modular braided tensor category, generalizing
thereby the well-known discriminant technique of quadratic form theory; cf.~\cite{Ho-shadow}. 

It should be mentioned that there are papers, which do not
distinguish between VOAs and the associated modular braided tensor category.
Usually, the authors make Assumption~\ref{hauptvermutung} given below and work 
in the setting of 3d-TQFTs. The example of the lattice VOAs
associated to the two positive definite even unimodular lattices $E_8^2$ and 
$D_{16}^+$ in dimension $16$~\cite{kneser} shows that there are
different VOAs which can lead to the same 3d-TQFT and which have the
same central charge.

\medskip

The paper is organized as follows.
In the second section, some basic definitions and results from the theory 
of integral quadratic forms are recalled. 
In section three, we define the genus of a vertex operator algebra
and try to formulate theorems or at least conjectures analogous to the ones
known for quadratic forms. 
In the final section, we reformulate the extension problem
for vertex operator algebras in terms of the associated modular 
braided tensor category as already outlined
in my paper~\cite{Ho-shadow}.
This extension problem can be considered as a kind of problem in coding 
theory over 3d-TQFTs.
The solution for VOAs defining an abelian intertwining algebra is given. 
Also, complete results for self-dual framed VOAs
or, equivalently, self-dual codes over the 3d-TQFT associated to the 
Ising model, are presented.

\smallskip

There are many open problems about the structure of vertex algebras
compared to what is known for quadratic forms.
But the theory of quadratic forms is more than two hundred years old. 
It seems to me that from the fifteen years old theory of vertex algebras
we can see only the tip of the iceberg of such a structure theory.

\medskip

I like to thank the organizers of the conference for inviting me to
the Fields Institute and for the possibility to present these ideas. 
I am grateful to Markus Rosellen for his comments on an early version
of the paper and for discussions about abelian intertwining algebras
as well as to Geoffrey Mason for reminding me of some references.
I also like to thank Jim Lepowsky and the referee for their many 
valuable suggestions. Finally, my thanks goes to Supun Koranasophonpun
for proofreading the manuscript.


\section{Genera of integral quadratic forms} 

We collect some results from the arithmetic of integral quadratic forms.
The theory of genera for binary forms was developed by Gau\ss\ in~\cite{gauss}. 
Instead of the more classical notion of quadratic forms in coordinates, 
we use the language of lattices.

\smallskip

We denote in this section by $\Z/N\Z$ the additive group of the integers modulo~$N$ and by
$\Z_p$ the ring of $p$-adic integers.
Recall that a map $q$ from an abelian group $B$ into
an additively written group $G$ is called a quadratic form 
if the expression $q(x+y)-q(x)-q(y)$ is additive
in $x$ and in $y$, i.e., bilinear.

\smallskip
{\bf Definition 2.A} \ 
A {\em finite quadratic space \/} $(A,q)$ is a finite abelian
group $A$ together with a quadratic form $q: A\longrightarrow \Q/2\Z$,
such that the induced bilinear form $b: A\times A \longrightarrow \Q/\Z$,
$b(x,y)=\frac{1}{2}(q(x+y)-q(x)-q(y))$ is nondegenerate.
Note that different quadratic forms can induce the same bilinear form.
\smallskip

\begin{theorem}[Structure theorem for finite quadratic spaces; cf.~\cite{Ni-genus}, Sec.~8]
A finite quadratic space can be decomposed uniquely as a direct sum:
$$(A,q)=\bigoplus_{p} (A_p,q_p),$$
where the $A_p$ are the $p$-Sylow subgroups of $A$ and $q_p$ is the 
quadratic form $q$ restricted to $A_p$.
Furthermore, each $(A_p,q_p)$ is the direct sum of quadratic spaces as follows:

For primes $p>2$, it is the sum of spaces $(\Z/{p^k}\Z,q_{\theta})$, with $k\geq 1$,
${\theta}\in \{\pm 1\}$ and 
quadratic form $q_{\theta}:\Z/{p^k}\Z=\langle x\rangle \longrightarrow \Q/2\Z$,
$q(x)=2\, c_{\theta}/p^k$, where $c_{\theta}$ is an integer with Jacobi symbol 
$\big(\frac{ c_{\theta}}{p}\big)=\theta$.
The only relations in the abelian semigroup  generated by the $(\Z/{p^k}\Z,q_{\theta})$
for fixed~$p$ arise from
$(\Z/{p^k}\Z,q_{\theta})\oplus(\Z/{p^k}\Z,q_{\theta})\cong
    (\Z/{p^k}\Z,q_{\theta'})\oplus(\Z/{p^k}\Z,q_{\theta'})$,
with ${\theta}$, ${\theta}'\in \{\pm 1\}$.

For $p=2$, one takes $\theta\in\Z_2^*/(\Z_2^*)^2$, there
are also two other kind of generators and further relations; cf.~\cite{Ni-genus},
Prop.~1.8.1 and~1.8.2.
\end{theorem}

\smallskip

{\bf Definition 2.B} \ An even lattice is a free $\Z$-module of 
finite rank together with a nondegenerate symmetric bilinear map
$(\,.\,,\,.\,): L\times L \longrightarrow {\bf Z}$,
such that $(x,x)\in 2{\bf Z}$ for $x\in L$.

\smallskip

The map $(\,.\,,\,.\,)$ can be linearly extended to $L\otimes R$
for any ring $R\supset {\bf Z}$.
Every lattice $L$ defines a quadratic space $(A,q)$ by letting
$A=L^*\!/L$, where $L^*=\{x\in L\otimes\Q\mid (x,y)\in \Z \ \hbox{for all } y\in L\}$
is the dual lattice, and $q:L^*\!/L\longrightarrow \Q/2\Z$ is
the quadratic form $x \pmod{L}\mapsto q(x)=(x,x) \pmod{2\Z}$,
called the discriminant form.
The {\it signature ${\rm sign}(L)$} of $L$ is the pair $(n_+,n_-)$,
where $n_+$ respectively~$n_-$ are the maximal dimension of a positive 
respectively~negative definite subspace of $L\otimes\R$. By Sylvester's Law of
inertia, $n_++n_-={\rm rk}(L)$, the rank of $L$.

\smallskip
{\bf Definition 2.C} \ The {\it genus of $L$\/} is the collection of the 
local lattices $L\otimes \Z_p$, $p$ a prime number, including
$L\otimes \R$ for $p=\infty$. 
\smallskip

It follows easily that the genus of $L$ determines $(L^*\!/L,q)$ and ${\rm sign}(L)$.
Also, the converse is true:

\begin{theorem}[cf.~\cite{Ni-genus}, Cor.~1.9.4]
The discriminant form $(L^*\!/L,q)$ and the signature ${\rm sign}(L)$ 
determine the genus of $L$.
\end{theorem}

Not every pair consisting of a finite quadratic space $(A,q)$
and a pair $(n_+,n_-)$ of nonnegative integers can be realized as the genus of 
a lattice.
First, there is a condition modulo $8$:
\begin{theorem}[Milgram]
For a lattice $L$ with discriminant form $(L^*\!/L,q)$ and signature $(n_+,n_-)$, one has
\begin{eqnarray}\label{milgram}
\frac{1}{\sqrt{|L^*\!/L|}}\sum_{x\in L^*\!/L} e^{\pi i q(x)} & =  & e^{2\pi i\, (n_+-n_-)/8}.
\end{eqnarray}
\end{theorem}

The next result shows that for ranks large enough, this is the only 
condition that the quadratic space and the signature have to satisfy.

\begin{theorem}
An even lattice $L$ with discriminant form $(A,q)$ and signature $(n_+,n_-)$ exists
if condition~(\ref{milgram}) is satisfied and $n_++n_-> {\rm rank}(A)$.
\end{theorem}

Note that $n_++n_-\geq {\rm rank}(A)$ is clearly a necessary condition.
For $n_++n_-={\rm rank}(A)$, one also has a precise but more complicated
condition for the existence of a lattice; see~\cite{Ni-genus}, Th.~1.10.1.

\begin{theorem}[Finiteness of the class number]
The set 
$${\rm gen}(L):=\{\hbox{\rm Isomorphy classes of lattices with the same genus as $L$}\}$$ 
is finite.
\end{theorem}

The proof uses the result that there are only finitely many lattices $L$ of fixed
discriminant $|L^*\!/L|$; see for example~\cite{Cassels}, Ch.~9, Th.~1.1.

Much more precise information about the number of lattices in a genus is known.
We restrict ourselves to positive definite lattices, i.e., lattices~$L$ of 
rank \hbox{$n=n_+$.} 

\begin{theorem}[Mass formula]
$$\sum_{M\in {\rm gen}(L)}\frac{1}{|{\rm Aut}(M)|}=2\,\cdot\, \alpha_{\infty}\times 
\prod_{p\,=\,2,\,3,\,5,\,\ldots}\alpha_p,$$
where $\alpha_{\infty}=2\pi^{-n(n+1)/4}\big(\prod_{j=1}^n\Gamma(j/2)\big)\cdot
|L^*\!/L|^{(n+1)/2}$ and
 the $\alpha_p$'s can be computed explicitly and depend only on $n$ and the $p$-component
$((L^*\!/L)_p,q_p)$.
\end{theorem}

This result was obtained independently by H.~J.~S.~Smith~\cite{Smith-quad}
and H.~Min\-kow\-ski~\cite{min}. 

The number of embeddings of lattices $K$ of rank $g\leq n$
into $L$ for all $K$ can be encoded
in the Siegel Theta series $\Theta^{(g)}_L(Z)$ of genus~$g$, a function on the 
Siegel upper half plane $\H_g=\{Z\in {\rm Mat}(g,\C)\mid Z=Z^t,\ {\rm Im}(Z)>0\}$ of genus $g$:
$$\Theta^{(g)}_L(Z)=\sum_{v_1,\,\ldots,\,v_g\in L}e^{2\pi i\,{\rm tr} 
    ((v_1,\ldots,v_g)(v_1,\ldots,v_g)^tZ) }.$$
It is a modular form of weight $n/2$ for a congruence subgroup $\Gamma$ of 
${\rm Sp}(g,\Z)$ acting on $\H_g$.

\begin{theorem}[Siegel]
$$\left( \sum_{M\in {\rm gen}(L)}\frac{\Theta_M^{(g)}(Z)}{|{\rm Aut}(M)|}\right)\,\cdot\,
\left(\sum_{M\in {\rm gen}(L)}\frac{1}{|{\rm Aut}(M)|}\right)^{-1}=
E^{(g)}(Z),$$
where $E^{(g)}(Z)$ is an explicitly given {\it Eisenstein series\/} for $\Gamma$,
depending only on the genus of~$L$.
\end{theorem}

This theorem was first proven by C.~L.~Siegel~\cite{sie}  and reformulated in an adelic
picture by A.~Weil, cf.~\cite{weil}, see also~\cite{Kudla-seesaw} for a general setting.

\medskip

There is a slightly finer invariant for lattices, 
namely M.~Eichler's {\it spinor genus\/},
cf.~\cite{Cassels}.
The set of lattices in a genus can be decomposed into $2^r$ spinor genera, 
also very often $r=0$.
The sum $\sum_M\frac{1}{|{\rm Aut}(M)|}$ over all lattices $M$ in the same spinor
genus as a  positive definite lattice $L$ is for all spinor genera in the 
genus of $L$ the same if ${\rm rank}(L)\geq 3$.
 
For indefinite lattices of rank at least $3$, there is only one lattice in 
each spinor genus. This is one way to see that there is only one
even unimodular (i.e., of discriminant $1$) lattice $I\!I_{r,s}$ of 
signature $(r,s)$ if ${\rm min}(r,s)\geq 1$. Furthermore, one can deduce that 
two lattices $L$ and $M$ are in the same genus if and only if $L\oplus I\!I_{1,1}$
and $M\oplus I\!I_{1,1}$ are isomorphic.

There is also the concept of {\it rational equivalence\/}. Two lattices $L$ and $M$
are called rational equivalent if the rational quadratic spaces 
$L\otimes \Q$ and $M\otimes \Q$
are isomorphic. An equivalent more geometric formulation is to say that both
$L$ and~$M$ have an isomorphic sublattice $K$ of finite index. 
{}From the Hasse principle, it follows that $L$ and $M$ are rational equivalent
if and only if all the local quadratic spaces $L\otimes \Q_p$ and $M\otimes \Q_p$, for $p$ 
a prime including $p=\infty$, are isomorphic. 
In particular, lattices in the same genus are rational equivalent.

Finally, we remark that the results explained in this section hold
similarly for non-even integral lattices.


\section{Genera of vertex operator algebras}

In this section, we explain how the concepts from the arithmetic of
quadratic forms as explained in the last section can, at least partially, 
be generalized to vertex operator algebras.

\smallskip

A good starting point for the description of three-dimensional topological
quantum field theories are modular braided tensor categories~\cite{Turaev}.
With the help of the Kirby calculus~\cite{Kirby}, one can use them to define 
invariants of $3$-manifolds, cf.~\cite{ReTu}.
We will give the basic definitions, for details see~\cite{Turaev}.

\smallskip
{\bf Definition 3.A}\  A {\it monoidal\/} category is a category 
${\mathcal{C}}$ together with a functorial associative tensor product 
$\otimes: \mathcal{C}\times \mathcal{C}\longrightarrow \mathcal{C}$
and a neutral object ${\bf 1}$ for $\otimes$. We assume that the tensor category
is {\it strict\/}, so for all objects $U$, $V$, $W$ of $\mathcal{C}$
the products $(U\otimes V)\otimes W$ and $U\otimes(V\otimes W)$ are identical
and not just isomorphic objects.

A {\it ribbon category\/} is a monoidal category together with:  
\begin{enumerate}
\item[(i)] Functorial isomorphisms $R_{U,V}:U\otimes V\longrightarrow V\otimes U$
(the {\it braiding\/}) such that 
$$R_{U,V\otimes W}=({\rm id}_V \otimes R_{U,W})(R_{U,V}\otimes {\rm id}_W)\ \hbox{and}\   
R_{U\otimes V, W}=(R_{U,W}\otimes {\rm id}_V)({\rm id}_U \otimes R_{V,W}).$$
\item[(ii)] Functorial isomorphisms $\theta_V: V\longrightarrow V$ (the {\it twist\/}) 
such that
$$\theta_{U\otimes V}=R_{V,U}(\theta_V\otimes\theta_U)R_{U,V}.$$
\item[(iii)] A triple $(*,d,b)$, which associates to any $V$ a {\it dual object\/} $V^*$ and 
morphisms $d: V^*\otimes V\longrightarrow {\bf 1}$ (the {\it evaluation\/})
and  $b: {\bf 1}\longrightarrow V\otimes V^*$ (the {\it coevaluation\/}), such that 
$$(d_V\otimes {\rm id}_{V^*})({\rm id}_{V^*}\otimes b_V)={\rm id}_{V^*} \  \hbox{and} \ \; 
({\rm id}_{V^*}\otimes d_V)( b_V\otimes {\rm id}_{V})= {\rm id}_V,$$
$$(\theta_V\otimes{\rm id}_{V^*})b_V=({\rm id}_{V}\otimes \theta_{V^*})b_V.  $$
\end{enumerate}

These axioms can be visualized by labeled bands in ${\bf R}^3$. 
In a ribbon category, one can define for any endomorphism $f:V\longrightarrow V$ 
the {\it trace\/} of $f$ as an element of ${\rm End}({\bf 1})$:
$${\rm tr}(f)=d_V R_{V,V^*} (\theta_V f \otimes {\rm id}_{V^*}) b_V.$$
The {\it quantum dimension\/} ${\rm dim}(V)$ of $V$ is the 
trace of the identity morphism of~$V$.

A {\it modular braided tensor category\/} (over $\C$), or 
{\it modular category\/} for short, is a ribbon category which has also the structure 
of an abelian category over $\C$ such that:
\begin{enumerate}
\item[(i)] the tensor product is $\C$-linear and ${\rm End}({\bf 1})=\C$;
\item[(ii)] every object is a direct sum of a finite set of simple objects;
\item[(iii)] the isomorphism classes of simple objects form a finite set $\{V_i\}_{i\in I}$
and an object $V$ is simple precisely if ${\rm End}(V)=\C$;
\item[(iv)] the matrix $({\rm tr}(R_{V_j,V_i}R_{V_i,V_j}))_{i,j \in I}$
is invertible.
\end{enumerate}

We call the number $D=\sum_{i\in I}({\rm dim}(V_i))^2$ the {\it discriminant\/}. 
In~\cite{Turaev}, a square root $\sqrt{D}$ is called a rank.
The $S$-matrix is defined as the matrix 
$S^{\rm TOP}=\frac{1}{\sqrt{D}}({\rm tr}(R_{V_j,V_i}R_{V_i,V_j}))_{i,j \in I}$.
We also need $\gamma=\left(\sum_{i\in I}\theta_i\,({\rm dim}(V_i))^2\big/
\sqrt{D}\right)^{1/3}$.
Here, the numbers $\{\theta_i\}_{i\in I}$ are defined by 
$\theta_{V_i}=\theta_i\cdot{\rm id}_{V_i}$.
The $T$-matrix is the diagonal matrix $T^{\rm TOP}$ with the 
$\gamma^{-1}\cdot\theta_i$, $i\in I$, on the diagonal.
The $S$- and $T$-matrix satisfy the relations $(ST)^3=S^2$, $S^2T=TS^2$ and $S^4=1$.
These are the same relations as for the generators 
$s=\bigl({{\phantom{-}0 \ 1}\atop{-1\ 0}}\bigl)$ 
and $t=\left({1 \ 1}\atop{0\ 1}\right)$ of ${\rm SL}_2(\Z)$,
the mapping class group of a genus~$1$ surface.
Therefore the $S$- and $T$-matrix define a complex representation~$\rho$ of ${\rm SL}_2(\Z)$ of 
dimension $|I|$.
If we omit the factor $\gamma^{-1}$ in the definition of the $T$-matrix, we obtain a
more natural projective representation of ${\rm SL}_2(\Z)$,
cf.~\cite{Turaev}, II.3.9. 
In the context of modular categories, one can prove that the $\theta_i$
and $\gamma$ are roots of unity; cf.~\cite{Va-rational,AnMo,DML2000,BaKi}.

Finally, there is the notation of a {\it unitary modular braided tensor category\/} 
where one has natural maps
$$\overline{\phantom{X}}:{\rm Hom}(V,W)\longrightarrow {\rm Hom}(W,V)$$
identical to complex conjugation for $V=W=V_i$, $i\in I$, such that
${\rm tr}(f\,\overline{f})\geq 0 $ for all $f\in{\rm Hom}(V,W)$.  
\smallskip 

Two modular categories are called equivalent, if there is a functorial
isomorphism between them, carrying over all structures. 

One can also consider nonstrict monoidal categories and define modular
categories by using them. Every monoidal category is equivalent to a strict one
by MacLane's coherence theorem and this is also true for modular braided tensor 
categories, cf.~\cite{Turaev},~XI.1.4.
The categories arising from conformal field theory are in general not 
strict. 

\smallskip

We formulate two versions of the {\it Verlinde formula\/} which are theorems in
the context of modular categories. 
First, let $(N_{ij}^k)^{\rm TOP}={\rm dim}\,{\rm Hom}(V_i\otimes V_j, V_k) $ be
the structure constants of the {\it fusion algebra\/} defined on the free vector 
space $\C[I]$ by the product $i\times j=\sum_{k\in I}(N_{ij}^k)^{\rm TOP} \,k.$ Then one has
$$(N_{ij}^k)^{\rm TOP}=\sum_{l\in I}\frac{S_{il}S_{jl}S_{kl^*}}{S_{0l}},$$
where $i^*\in I$ is the label defined by $V_{i^*}\cong V_i^*$ and $0\in I$ is the label for
$V_0\cong {\bf 1}$. Secondly, let $N^{\rm TOP}_{i_1,\ldots,i_n}(g)$
be the dimension of the vector space associated by the
$3d$-TQFT to a genus~$g$ surface decorated with labels $i_1$, $\ldots$, $i_n$.
Then one has
$$ N^{\rm TOP}_{i_1,\ldots,i_n}(g)=D^{g-1}\sum_{j\in I}\left(
({\rm dim}(V_j))^{2-2g-n}\prod_{s=1}^nS_{i_s,j}\right) .$$

We mention three important classes of examples of modular categories:

1) The modular category $\mathcal{C}(A,q)$ associated to a finite quadratic space $(A,q)$,
see~\cite{Turaev}, I.1.7.2 and~II.1.7.2. 
We have $q:A\longrightarrow \Q/2\Z\subset U(1)\subset\C^*$, $q(x)=\varphi(x)c(x,x)$,
where the maps $c:A\times A\longrightarrow \C^*$, $\varphi: {}_2A\longrightarrow \C^*$
(${}_2A$ denotes the $2$-torsion subgroup of $A$)
used in~\cite{Turaev} are up to equivalence identical to the map~$q$, 
see~\cite{FrKe}, Sec.~7.5.
Since the fusion algebra of $\mathcal{C}(A,q)$ is just the group ring $\C[A]$
of the abelian group $A$, we call such a modular category abelian. 
Any modular category with such an abelian group ring $\C[A]$
as fusion algebra arises in such a way, see~\cite{FrKe}, Ch.~7.

2) The Chern-Simons theory associated to a simple Lie group $G$ and 
positive integral level~$k$.
This modular category can be constructed using quantum groups, 
cf.~\cite{Drinfeld,ReTu}.

3) Discrete Chern-Simons theories constructed from a finite group $G$ and 
a class~\hbox{$c\in H^3(G,{\C})$}~\cite{DiWi,FrQu}.

\begin{prob}
Find a good description of the set of equivalence classes
of modular braided tensor categories over ${\bf C}$.
\end{prob}

For two such categories $\mathcal{C}$ and $\mathcal{D}$, it is easy to define 
a product category $\mathcal{C}\times \mathcal{D}$, such 
that its simple objects are the products of
the simple objects of the factors.
One can decompose a category into indecomposable pieces $\mathcal{C}_s$
under the product: $\mathcal{C}=\mathcal{C}_1\times\mathcal{C}_2\times\ldots \times\mathcal{C}_r$.
An additional problem compared to the situation for quadratic spaces is that 
it is not possible to decompose $\mathcal{C}$ into 
a product of local pieces $\mathcal{C}_p$
for different primes $p$ such that,
for example, the numbers $\theta_i$ for $\mathcal{C}_p$ are $p^e$-th root of unity.

A conjecture for modular braided tensor categories (also usually formulated 
for rational VOAs) says that the above introduced representation
$\rho:{\rm SL}_2(\Z)\longrightarrow {\rm End}(\C^{|I|})$ has a congruence subgroup as
kernel or, more precisely, $\ker(\rho)\supset \Gamma(N)$ where
$N$ is the smallest natural number such that $T^{N}={\bf 1}$.
Away from the prime $p=2$, this has been proven by Coste and Gannon~\cite{Gannon-modular}.
P.~Bantay~\cite{Ban-mod} has recently reduced the problem
to a property of permutation orbifolds. 
Assuming this conjecture, Eholzer~\cite{Eh-dr} has classified the  
indecomposable modular categories (or at least the fusion algebras together with 
the $T$-matrix and $\gamma$) with  $ |I|\leq 4 $ using the representation theory 
of ${\rm SL}_2(\Z/N\Z)$.
A lot of work has been done on the classification of fusion algebras and
also compatible tensor categories, cf.~for example~\cite{FrKe,Fu-fusion,MN2001}.

\medskip

For the purpose of this paper, the following definition of a 
restricted class of vertex algebras seems to be most adequate.

\smallskip
{\bf Definition 3.B}\ A {\it vertex operator algebra of
central charge $c$\/} is a graded complex 
vector space \hbox{$V=\bigoplus_{n\in {\bf Z}_{\geq 0}} V_n$} with ${\rm dim}\, V_n<\infty$ together with a 
linear map $Y(\,.\,,z):V\longrightarrow {\rm End}(V)[[z,z^{-1}]]$ 
({\it the state field correspondence})
such that for $v\in V_m$ the coefficients $v_n$ of $Y(v,z)=\sum_{n\in\Z} v_n z^{-n-1}$
are of degree $m-n-1$ and $v_nu=0$ for fixed $u\in V$ and $n$ large enough,
and there are two distinguished elements ${\bf 1}\in V_0\cong\C$ (the {\it vacuum\/}) and $\omega\in V_2$ 
(the {\it Virasoro element\/}) subject to the following axioms:
\begin{itemize}
\item[-] $Y({\bf 1},z)={\rm id}_V$ and for $v\in V$, one has $v_n{\bf 1}=0$ for $n\geq 0$, $v_{-1}{\bf 1}=v$.
\item[-] For  $u$, $v\in V$, there is a $N\in \Z_{\geq 0}$ with 
$(z-w)^N Y(u,z)Y(v,w)=\linebreak 
    (z-w)^N Y(v,w)Y(u,z)$ in ${\rm End}(V)[[z,z^{-1},w,w^{-1}]]$.
\item[-] The coefficients of 
$Y(\omega,z)=\sum_{n\in\Z} L_n\, z^{-n-2} $
define a representation of the {\it Virasoro algebra\/} 
of central charge $c$: $[L_m,L_n]=L_{m+n}+\frac{m^3-m}{12}\delta_{m+n,0}\cdot c\cdot  {\rm id}_V$
and for $v\in V$, one has $Y(L_{-1}v,z)=\frac{d}{dz} Y(v,z)$.
\end{itemize}
\smallskip

These axioms can be distilled from the general axioms of axiomatic 
quantum field theory applied to the case of Riemannian surfaces as
the space time. To obtain a more geometric formulation
as indicated in the introduction from a VOA, further conditions 
have to be satisfied.
For the terms we use below and for details, we refer to the literature: 

There is the notion of modules and intertwining spaces 
$\left(M_3 \atop M_1,\ M_2\right) $ for a VOA, see~\cite{FHL}.
A VOA is usually called {\it rational\/}, if there are only finitely 
many nonisomorphic irreducible modules and any module can be decomposed into
a finite direct sum of irreducible ones. 
(The first condition can be deduced from the second, cf.~\cite{DML98}.)
The  {\it conformal weight\/} $h$ of a module 
is its smallest $L_0$-eigenvalue. Let $\{M_j\}_{j\in I}$ be a complete set 
of nonisomorphic irreducible modules.
The $T$-matrix is defined as the diagonal matrix $T^{\rm VOA}$ with the 
numbers $e^{- 2\pi i c/24}\cdot e^{2\pi i h_j}$, $j\in I$, on the diagonal.
Zhu has shown~\cite{Zhu-dr} that under some conditions on the rational VOA, 
there is a representation of ${\rm SL}_2({\bf Z})$ on some space of genus~$1$
correlation functions. The generator $s$ of ${\rm SL}_2({\bf Z})$
defines the matrix $S^{\rm VOA}$.
Using the tensor product theory for modules of VOAs, Huang and Lepowsky 
associated a (nonstrict) braided tensor category 
to a suitable class of VOAs, see~\cite{HuLe-tensor} and the
references in~\cite{HuLe-affine}.  Finally, there is a definition for  
the dimension $N^{\rm VOA}_{i_1,\ldots,i_n}(g)$  of the space of vacua
on the genus~$g$ surface with labels $i_1$, $\ldots$, $i_n$ by Zhu~\cite{Zhu-riem}.

I assume further that the VOA is {\it unitary\/}, this means the
VOA can be defined over the real numbers and the natural invariant symmetric form
on it is positive definite, although this will exclude many interesting 
examples. Also, the VOA should be {\it simple\/}, i.e., $V$ is irreducible
as a module over itself.

We like to use the above mentioned results and some other properties believed
to be true for many interesting VOAs. We therefore make the following
assumption on the VOAs considered:

\begin{assumption}\label{hauptvermutung}
The VOA $V$ satisfies the conditions needed in the construction of 
the braided tensor category on the category of modules as in~\cite{HuLe-tensor} 
and also satisfies the conditions used in~\cite{Zhu-dr,Zhu-riem}.

Furthermore, the braided tensor category of $V$-modules has the structure of a  
modular braided tensor category 
$\mathcal{C}(V)$ such that 
$$\begin{array}{lrcl}
\hbox{(i)}  \qquad \qquad \qquad & e^{2\pi i c/24}\cdot T^{\rm VOA} & = & \gamma \cdot T^{\rm TOP}, \\
\hbox{(ii)} \qquad \qquad \qquad  & S^{\rm VOA} & = & S^{\rm TOP}, \\
\hbox{(iii)}  \qquad \qquad \qquad  & \dim \left({M_k}\atop{M_i,\,M_j}\right) & = & 
     \dim {\rm Hom}(V_i\otimes V_j,V_k),\ \ \hbox{for $i$, $j$, $k\in I$,} \\
\hbox{(iv)} \qquad \qquad \qquad  & N^{\rm VOA}_{i_1,\ldots,i_n}(g) & = & N^{\rm TOP}_{i_1,\ldots,i_n}(g),
\ \ \hbox{for $i_1$, $\ldots$, $i_n\in I$, $g\geq 0$.}
\end{array}$$
\end{assumption}

One problem is certainly to find simple sufficient conditions on the VOA so that the
assumption becomes a theorem.
Sometimes other conditions in the definition of rationality are added.
Assuming that Assumption~\ref{hauptvermutung}
is true would directly imply that the central charge and 
conformal weights are rational numbers.

Property (iv) follows from (iii) if one can prove the sewing relations
$$\sum_{k\in I}  N^{\rm VOA}_{i_1,\ldots,i_n,k}(g)\cdot N^{\rm VOA}_{k^*,j_1,\ldots,j_m}(h)
= N^{\rm VOA}_{i_1,\ldots,i_n,j_1,\ldots,j_m}(g+h)$$
and 
$$\sum_{k\in I}  N^{\rm VOA}_{i_1,\ldots,i_n,k,k^*}(g)=
                           N^{\rm VOA}_{i_1,\ldots,i_n}(g+1).$$

For important examples of VOAs, at least parts of the assumption are known to be true;
for details we refer to the references:

1) VOAs defining an abelian intertwining algebra; cf.~\cite{DoLe,Ro-aia}.
     
2) Wess-Zumino-Witten models, i.e., the VOAs defined on the
highest weight representations of level $k$ for an affine Kac-Moody 
Lie algebra $\widetilde{{\rm Lie}(G)}$; cf.~\cite{TUY,KaLu,FreZhu,Faltings,Finkel,HuLe-affine}.

3) Holomorphic orbifolds; cf.~\cite{DVVV,Hua,DoMa-moonshine}.

4) Minimal models; cf.~\cite{Hu-min,HuMi-minI,HuMi-minII}.

\smallskip

One could try to generalize the above mentioned results to rational vertex algebras.
In their definition one allows $\Z$-graded vector spaces $V$ and the 
pieces $V_n$ are not assumed to be finite-dimensional. Vertex algebras should be
considered as the analog of indefinite even lattices.
For a vertex algebra $V$ of central charge $c$, one can define the number $c_+$
as the largest central charge of a subVOA $W$ of $V$, $c_-$ as
$c-c_+$ and call the pair $(c_+,c_-)$ the signature of~$V$. 

\smallskip

For the rest of this section, we consider only VOAs for which  
Assumption~\ref{hauptvermutung} is satisfied. 

\smallskip
{\bf Definition 3.C} \  The {\it genus\/} of a rational VOA $V$ of 
central charge $c$ is defined as the pair consisting of the associated 
modular braided tensor category $\mathcal{C}(V)$ and the central charge~$c$.
\smallskip

It is not known much about the set of central charges which can occur.
{}From the discussion of modular categories, we know that $c$ has to be
rational. The theory of unitary
highest weight representations of the Virasoro algebra shows  that
$c$ belongs to the minimal series $c=1-\frac{6}{n(n+1)}$, $n=3$, $4$, $\ldots$,
for $c<1$.

There is the following analogue of Milgrams theorem for modular braided
tensor categories.

\begin{theorem}[cf.~\cite{Turaev}]\label{thVOAmilgram}
For a rational VOA $V$ with associated modular category $\mathcal{C}(V)$ and 
central charge $c$, one has
\begin{eqnarray}\label{VOAmilgram}
\frac{1}{\sqrt{\sum_{i\in I}\dim(V_i)^2}}\sum_{j\in I}e^{2\pi i h_j}\dim(V_j)^2
=e^{2\pi i c/8}.
\end{eqnarray}
\end{theorem}

{\bf Proof:}
The result is an immediate consequence of Assumption~\ref{hauptvermutung} and
the existence of the ${\rm SL}_2({\bf Z})$-representations.
We have $(S^{\rm VOA} T^{\rm VOA})^3=(S^{\rm VOA})^2$ 
and $(S^{\rm TOP} T^{\rm TOP})^3=(S^{\rm TOP})^2$ because of the ${\rm SL}_2(\Z)$-representations.
Using properties (i) and (ii) in Assumption~\ref{hauptvermutung}, we get 
$e^{2\pi i c/8}=\gamma^3$.
\qed

\begin{prob}
Let $\mathcal{C}$ be a modular braided tensor category. Does there exist rational
VOAs $V$ of central charge $c$ with $\mathcal{C}(V)=\mathcal{C}$ if  
condition (\ref{VOAmilgram}) is satisfied and $c$ is large enough?
\end{prob}

What are the possible values of $c$ for abelian $\mathcal{C}\cong \mathcal{C}(A,q)$? 
For VOAs there are more central charges $c$ possible than one can realize by 
lattice VOAs.
The bound $c\geq {\rm rank}(A)$ does not hold in general:
The fixpoint subVOA $V_{\sqrt{2}E_8}^\tau$
of the lattice VOA $V_{\sqrt{2}E_8}$ 
for the involution $\tau$ lifted from the reflection at the origin of $\sqrt{2}E_8$
(cf.~\cite{Gr-baby}) has central charge~$8$, but 
$\mathcal{C}(V_{\sqrt{2}E_8}^\tau)\cong\mathcal{C}({\bf F}_2^{10},q)$, where
$q$ is isomorphic to $5$ copies of the hyperbolic plane.

\begin{conj}
The set 
$${\rm gen}(V):=\{\hbox{\rm Isomorphy classes of VOAs with the same 
genus as $V$}\}$$ 
is finite.
\end{conj}

{\it Example:\/} Let $\mathcal{C}$ be the trivial modular category, i.e., there is only one
simple object. The corresponding VOAs are called self-dual or holomorphic.
By Theorem~\ref{thVOAmilgram}, the central charge $c$ has to be divisible by $8$.
An interesting case is $c=24$. Schellekens has found a list of $71$ candidates
for the corresponding genus~\cite{schellekens1}.
For the $24$ lattice VOAs constructed from the genus of the positive definite even
unimodular lattices in dimension $24$, one has an existence and uniqueness theorem. 
$24-9$ other VOAs are completely constructed as $\Z_2$-orbifolds, cf.~\cite{DGM-twist,DGM}. 
One of them is the Moonshine module $V^{\natural}$, cf.~\cite{Bo-ur,FLM}. 
The fact that $V^{\natural}$ belongs to 
the genus was proven in~\cite{Do}. The uniqueness conjecture~\cite{FLM} that $V^{\natural}$
is the only VOA $V$ in the genus with $\dim V_1=0$ is still open.
For the VOA with $V_{A_{1,2}}^{\otimes 16}$ as the affine Kac-Moody subVOA, 
we can show uniqueness, cf.~\cite{GH-stabil} and next section. 

Conformal field theory suggests that for $g\geq 1$ 
it may be possible to define a (possibly vector valued) genus 
$g$ partition function $\chi_V^{(g)}$ on the
genus $g$ Teichm\"uller space $\mathcal{T}_g$ with automorphic
properties for the mapping class group of a genus~$g$ surface. 
For $g=1$ this follows from  Zhu's work~\cite{Zhu-dr}. 

\begin{prob}
Are there good Mass formulas for the genus of a VOA of the form 
$$\sum_{W\in {\rm gen}(V)}\frac{1}{|{\rm Aut}(W)|} \quad \hbox{and}
\sum_{W\in {\rm gen}(V)}\frac{\chi_W^{(g)}(Z)}{|{\rm Aut}(W)|}\ ?$$
\end{prob}\nopagebreak

One problem is to find a good definition for $|{\rm Aut}(V)|$ if 
${\rm Aut}(V)$ is infinite; one could try to take the Tamagawa measure 
if ${\rm Aut}(V)$ is a reductive Lie group defined over ${\bf Q}$ by using
VOAs defined over $\Q$.
For a kind of relative Mass formula of this type involving 
a nonabelian extension problem, see Theorem~\ref{relmass} in the next section.

\smallskip

The definition of the genus of a VOA  was modeled along the notation of
the genus of an even lattice $L$ such that the following diagram is 
commutative
$$\begin{array}{ccc}
L & \longrightarrow & V_L \\
\downarrow & & \downarrow \\
{\rm genus}(L)=((L^*\!/L,q),{\rm rk}(L) )& \longrightarrow & {\rm genus}(V)=(\mathcal{C}(V_L),c(V_L)) .
\end{array}$$
The upper horizontal arrow is the lattice construction of VOAs~\cite{Bo-ur},
the left down arrow gives the genus of a lattice  via the discriminant form,
the right down arrow  comes from the abelian intertwining algebra structure on the
modules of $V_L$~\cite{DoLe} and the lower horizontal arrow is the map associating the
abelian modular category to a finite quadratic space as explained 
as an example above.

Regarding genera of vertex algebras, a structure theory seems to be 
more complicated than for indefinite lattices:  
The vertex algebras $V^{\natural}\otimes V_{I\!I_{1,1}}$ and 
$V_{\Lambda}\otimes V_{I\!I_{1,1}}\cong V_{I\!I_{25,1}}$ 
(here, $\Lambda$ denotes the Leech lattice) are not isomorphic
since both define in a uniform way different generalized Kac-Moody Lie algebras,
namely the Monster~\cite{Bo-Lie} and the fake Monster Lie algebra~\cite{Bo-fake}.
(They are different
because all Cartan subalgebras of a generalized Kac-Moody Lie algebra
are conjugate by inner automorphisms as shown by U.~Ray~\cite{Ray-cartan}.)
But both should be self-dual vertex algebras of signature $(25,1)$.
Another such vertex algebra is studied in~\cite{HS}.

A natural generalization of the notion of {\it rational equivalence\/} from lattices 
to VOAs would be as follows: 
We take the smallest equivalence relation such that 
two rational VOAs $V$ and $W$ of central charge $c$ are equivalent 
if both have isomorphic rational subVOAs with the same Virasoro element. 
This means that
there exist VOAs $U_0=V$, $U_1$, $\ldots$, $U_r=W$ and VOAs $S_0$, $\ldots$, $S_{r-1}$
such that $S_i$ is a common subVOA of $U_i$ and $U_{i+1}$ for $i=0$, $\ldots$, $r-1$.
Just to say there exists a common isomorphic subVOA wouldn't give an equivalence relation:
take for $V$ and $W$ the fixpoint VOAs $V=U^G$ and $W=U^{G'}$ where $U$ is the
$A_1\cong\sqrt{2}\Z$ lattice VOA and $G$ and $G'$ are finite subgroups of 
${\rm SU}(2)$ which can't be conjugated into the same larger finite subgroup of ${\rm SU}(2)$
(like the binary icosahedral and the binary tetrahedral group).

Also, one can generalize the discussion of this section to 
vertex operator super algebras (SVOAs). 
The central charge $c$ of a self-dual SVOA has to be half-integral~\cite{Ho-dr},
Th.~2.2.2.
For an interesting example with $c=23\frac{1}{2}$ and 
the classification for small $c$ see also~\cite{Ho-dr}.


\section{Extension problems}

In this section, we explain how extension problems for
suitable VOAs can be formulated in terms of the associated
modular braided tensor category and sometimes solved by this methods.

\medskip

For even lattices, all overlattices can easily be described; cf.~\cite{Ni-genus}:
\begin{theorem}[Lattice extensions]\label{latticeext}
\begin{itemize}
\item[(a)] The even overlattices $K$ of an even lattice $L$ are in one-to-one 
correspondence to the isotropic subspaces~$C$ of the associated quadratic space
$(L^*\!/L,q)$.
\item[(b)] The quadratic space for $K$ is given by 
the pair $(C^{\perp}/C ,(q\mid C^{\perp})/C)$.\nopagebreak
\item[(c)] Two such extensions give isomorphic lattices $K$ exactly if the corresponding 
subspaces can be transformed into each other by an isometry of $(L^*\!/L,q)$ induced 
from an automorphism of $L$. 
\end{itemize}
\end{theorem}

For integral lattices, one has to take subspaces $C$ on which the induced bilinear 
form \hbox{$b_L:L^*/L\times L^*/L\longrightarrow {\bf C}$} is trivial. 

\medskip

Our aim is to present similar results for extensions of VOAs.
Let $V$ be a rational VOA.
We call a VOA $(W,Y_W)$ an {\it extension\/} of $V$ if it
contains a subVOA isomorphic to $V$ and has the 
same vacuum and Virasoro element as~$V$. 
We have $W\cong\bigoplus_{i\in I}n_i\, M_i$ for a decomposition of $W$ 
into irreducible $V$-modules $M_i$,
where $\{M_i\}_{i\in I}$ is a complete set of irreducible modules 
of~$V$ and the nonnegative integers $n_i$ are the multiplicities of $M_i$ in $W$.
Also, we assume that we have fixed an embedding $\iota_W:V\longrightarrow W$.

We call two extensions $(W,Y_{W})$ and $(W',Y_{W'})$ isomorphic 
if there is an VOA-isomorphism $\sigma:W\longrightarrow W'$ and an 
automorphism $\tau$ of $V$ such that $\sigma\circ \iota_W=\iota_{W'}\circ \tau$.
For $(W,Y_{W})=(W',Y_{W'})$ and $\iota_W=\iota_{W'}$, 
we can define in this way
the {\it automorphism group\/} ${\rm Aut}_V(W)$ of the extension $(W,Y_{W})$.

Part (a) of Theorem~\ref{latticeext} implies that the overlattices
of a lattice can be determined from the structure of its
discriminant form. For VOAs we have as an analogous result:

\begin{theorem}[see~\cite{Ho-shadow}, p.~615]\label{voaext}
The VOA-extensions $W$ of a rational VOA $V$ satisfying 
Assumption~\ref{hauptvermutung} can be determined
completely in terms of the associated modular braided tensor category $\mathcal{C}(V)$. 
\end{theorem}

{\bf Proof:} Fix a $V$-module isomorphism $W\cong\bigoplus_{i\in I}n_i M_i$
extending $\iota_W$. The vertex 
operator $Y_W$ is an element of 
$$\bigoplus_{i,\,j,\,k\in I}{\rm Hom}({\bf C}^{n_i}\otimes 
{\bf C}^{n_j},{\bf C}^{n_k})\otimes \left({M_k}\atop{M_i,\ M_j}\right)_V.$$ 
For $(W,Y_W)$ to be a VOA-extension of $V$,
the modules $M_i$ with $n_i>0$ have to be of integral conformal weight and $Y_W$
has to satisfy the commutativity axiom. Both conditions are completely controlled 
by the twist and the braiding of the associated modular braided tensor category.
$\quad$ 
\qed

Similarly, one can see that the structure of the intertwining algebra of an 
extension is completely determined by $\mathcal{C}(V)$ and 
$Y_W$.~\footnote{During preparation of this manuscript, 
the paper~\cite{KiOs} by A.~Kirillov, Jr.~and V.~Ostrik became available.
The authors formulate the following precise analogues of part (a) and~(b) 
of Theorem~\ref{latticeext} (see Th.~6.2 of~\cite{KiOs}) using notation 
from~\cite{BaKi}:
(a) The extensions $W$ of~$V$ are in one-to-one correspondence to rigid
$\mathcal{C}(V)$-algebras $A$ on which the twist $\theta$ is trivial.
(b) The category of representations of $W$ can be identified with ${\rm Rep}^0 A$.}
We expect that any rational VOA $V$ has only finitely many nonisomorphic 
extensions $(W,Y_W)$.

In the case that $\mathcal{C}(V)$ is abelian, i.e., $\mathcal{C}(V)=\mathcal{C}(A,q)$ 
for a finite quadratic space $(A,q)$, one could expect in generalization
of Theorem~\ref{latticeext} (b) that the VOA-extensions are determined by 
the isotropic subspaces of $A$. This is indeed the case. 
Recall that a simple module $M_i$ is called a {\it simple current\/}
if for each simple module $M_j$ there is a another simple module $M_{j'}$
such that $M_i\times M_j=M_{j'}$ holds in the fusion algebra.

\begin{theorem}[Simple current extensions]\label{simplecurrent}
Let $V$ be a rational VOA which  has an abelian intertwining operator algebra
structure on the direct sum of the simple currents. 
Let $C$ be a subgroup of the abelian group $A\subset I$ of labels of the simple currents
for which the modules $M_c$, $c\in C$, have integral conformal weight. 
Then there exists 
a unique simple VOA-extension $(W,Y_W)$ of the form $W\cong\bigoplus_{c\in C}n_c\, M_c$,
$n_c>0$, and one has $n_c=1$.
\end{theorem}

{\bf Proof:}
{\it Existence:\/} Choose a non-zero vector
$\mathcal{Y}_{c,c'}\in \left(M_{c+c'}\atop M_c,\,M_{c'}\right)\cong{\bf C}$ for  
\hbox{$c$, $c'\in C$,} and define $Y_W=\bigoplus_{c,\,c'\in C}\mathcal{Y}_{c,c'}^{}$.
Recall~\cite{Mac-lane} that the first cohomology groups $H^i(K(C,2),\C^*)$ 
of the second  Eilenberg-MacLane space $K(C,2)$ (which is up to homotopy
the unique CW-complex with homotopy groups $\pi_2( K(C,2))=C$,
and $\pi_i( K(C,2))=0$ for $i\not=2$)
can explicitly described by the following cochain complex:
{\small
$$\begin{array}{ccccccccc}
K^5 &  \stackrel{\partial^4}{\longleftarrow} & K^4 & \stackrel{\partial^3}{\longleftarrow} & 
          K^3 &\stackrel{\partial^2}{\longleftarrow} & K^2 & \cdots \\
\\
\ldots & & F:C\times C\times C\rightarrow \C^* &  & f: C\times C\rightarrow \C^* 
&   & \varphi: C\rightarrow\C^* & \\
& & H:C\times C \rightarrow \C^*  & & & & & &
\end{array} $$ }
The coboundary maps $\partial^i$ which we need are defined by
 $\partial^2(\varphi)=\widehat\partial^2(\varphi)$ and 
$\partial^3(f)=(\widehat\partial^3(f),f(c,c')/f(c',c))$, where $\widehat\partial^i$ 
denotes the usual group cohomology coboundary operator.
Note that we can work with normalized cochains, i.e., cochains with 
$\varphi(0)=1$, $f(0,c)=f(c,0)=1$, respectively $F(0,c,c')=F(c,0,c')=F(c,c',0)=1$ and
$H(0,c)=H(c,0)=1$ for $c$, $c'\in C$.
The commutativity and associativity properties for abelian intertwining
algebras~\cite{DoLe} associate to $Y_W$ a cocycle $(F,H)\in K^4$.
The operator $Y_W$ defines a VOA structure on $W$ if $(F,H)=(1,1)$.
The cohomology group $H^4(K(C,2),\C^*)$ can be identified via $q'(x)=H(x,x)$
with the space of quadratic forms $q'$ on $C$, see~\cite{Mac-lane}, Th.~3.
In our situation, we have $q'=1$ since the modules $M_c$, $c\in C$, 
have integral conformal weight and $H(c,c)$ equals the exponential function applied
to $2\pi i$ times the conformal weight of $M_c$.
Therefore, there exists a $f\in K^3$ with $\partial^3(f)=(F,G)$. 
Replacing $\mathcal{Y}_{c,c'}$
by $f(c,c')^{-1}\mathcal{Y}_{c,c'}$, we see that $(W,Y_W)$ is a VOA.

{\it Uniqueness:\/} Assume $(W,Y_W)$ is an extension.
Similar as in the proof of Prop.~2.5 (4) of~\cite{DGH-virs}, one gets 
$m_c=m_c'$ for all $c$, $c'\in C$. But $m_0=1$ since $M_0$ must have multiplicity one.
One has $\mathcal{Y}_{c,c'}\not=0$ since by Proposition~11.9 of~\cite{DoLe} one has
$Y_W(u,z)v\not=0$ if $u$ and $v$ are not $0$.
We can replace $\mathcal{Y}_{c,c'}$ by  $f(c,c')\mathcal{Y}_{c,c'}$ 
where $f:C\times C\rightarrow \C^*$ is a $3$-cocycle, i.e., $\partial^3(f)=(1,1)$.
A $V$-module isomorphism of the space $W\cong\bigoplus_{c\in C}M_c$ defines  
for any irreducible $M_c$, $c\in C$, a scalar $\varphi(c)\in \C^*$, i.e., a two
cochain $\varphi\in K^2$, and $f$ is changed to $\partial^2(\varphi)f$.
But from the universal coefficient theorem, $H^3(K(C,2),\C^*)={\rm Ext}(C,\C^*)$ 
and ${\rm Ext}(C,\C^*)=0$ since the $\Z$-module $\C^*$ is divisible. So,
VOAs $(W,Y_W)$ for different choices of cocycles $f$ are isomorphic.
\qed

{\it Remarks:\/} In terms of the associated modular braided tensor category,
the simple currents define a subcategory $\mathcal{A}\subset \mathcal{C}(V)$ which 
is abelian, i.e., $\mathcal{A}\cong\mathcal{C}(A,q)$ for a finite quadratic 
space $(A,q)$ (we allow here that $q$ is degenerate,
i.e., $\mathcal{C}(A,q)$ has not to be modular):
The space of quadratic forms on $A$ can be identified as before
with the cohomology group $H^4(K(A,2),\C^*)$ and one has $q=H(x,x)$.
Furthermore, $q|C=1$ since the modules $M_c$, $c\in C$, 
have integral conformal weight, i.e, $C$ is an isotropic subspace of $(A,q)$.

\smallskip

The theorem completely solves the extension problem for VOAs $V$ for
which on the direct sum of a complete set of nonisomorphic irreducible modules
there is the structure of an abelian intertwining algebra. 
In this case, the modular braided tensor category $\mathcal{C}(W)$ 
of an extension $W$ are described as in Theorem~\ref{latticeext} for lattices.


The structure of the modular braided tensor category $\mathcal{C}(W)$ for 
general simple current extensions
was investigated in~\cite{FSS}, where a conjecture for the set of irreducible $W$-modules 
and the corresponding $S$-matrix was made. 
This conjecture was proven in~\cite{Mu}.

There is a natural action of the dual group $\hat C={\rm Hom}(C,\C)\leq{\rm Aut}_V(W)$ on~$W$.  

In the context of braided tensor categories, the
problem was also studied in \cite{FrKe}, section~7.6. 

\smallskip

Another special situation where an extension problem can be 
solved partially arises from finite group actions. On the VOA side, one has:

\begin{theorem}[Galois correspondence~\cite{DoMa-quantumgalois,HMT}]
Let $W$ be a VOA and $G\subset {\rm Aut}(W)$ be a finite subgroup of the  automorphism group.
Let $V=W^G$ be the fixpoint VOA. Then the extensions $U$ of $V$ which are 
contained in $W$
are in one-to-one correspondence to the subgroups of $G$, i.e., there is 
a subgroup $H\leq G$ such that $U=W^H$.
\end{theorem}

{}From the perspective of modular braided tensor categories, this situation was 
analyzed independently by~Brugui\`eres~\cite{Bru-cat} and M\"uger~\cite{Mu}.
Both give a purely topological characterization of the
above situation using a characterization result for neutral Tannakian categories obtained
by Doplicher and Roberts~\cite{DoRo} and Deligne~\cite{Deligne}, respectively. 
Both authors describe how to construct the modular braided tensor category
$\mathcal{C}(W)$ from $\mathcal{C}(V)$ and the representation subcategory for the group $G$. 
M\"uger also proves the Galois correspondence, see~\cite{Mu}, Sec.~4.2.

\smallskip

As a final example, we completely solve the extension problem for self-dual
framed vertex operator algebras (FVOAs). 
FVOAs were introduced in~\cite{DGH-virs} and are VOAs which have a
tensor product $L_{1/2}^{\otimes r}(0)$ of $r$ central charge $1/2$ rational Virasoro
VOAs as subVOA with the same Virasoro element, i.e., they are the extensions of 
$L_{1/2}^{\otimes r}(0)$. Recall that $L_{1/2}(0)$ has three irreducible modules
$L_{1/2}(0)$, $L_{1/2}(1)$ and $L_{1/2}(2)$ of conformal weight $0$, $\frac{1}{2}$
and $\frac{1}{16}$ and $L_{1/2}(1)$ is a simple current.
{}From the decomposition
$$V=\bigoplus_{h_1,\,\ldots,\,h_r\in\{0,1,2\}} n_{h_1,\,\ldots,\,h_r}\,
 L_{1/2}(h_1)\otimes\cdots \otimes L_{1/2}(h_r),$$
one gets the two codes 
$$\mathcal{C}=\{c\in\{0,1\}^r\mid n_{c_1,\,\ldots,\,c_r}\not=0\}\ \ \hbox{and} \ \ 
\mathcal{D}=\{d\in\{0,1\}^r\mid \sum_{ h_1,\ldots,\, h_r \atop h_i\,=\,2\
 \hbox{\scriptsize if}\ d_i\,=\,1} n_{h_1,\,\ldots,\,h_r}\not=0\}$$
and it was shown in~\cite{DGH-virs} that this codes
are linear, $\mathcal{D}\subset \mathcal{C}^{\perp}$ ($C^{\perp}$ denotes the code orthogonal 
to $C$), $\mathcal{C}$ is even,
and all weights of vectors in $\mathcal{D}$ are divisible by $8$.

Theorem~\ref{simplecurrent} says that there is a unique FVOA of the form 
$V_\mathcal{C}=\bigoplus_{c\in \mathcal{C} }L_{1/2}(c_1)\otimes\cdots\otimes L_{1/2}(c_r)$
for every linear even code $\mathcal{C}$, a result proven in~\cite{DGH-virs} 
(Prop.~2.16) and also obtained in~\cite{Mi-code}. As mentioned above, 
the associated modular category for $V_\mathcal{C}$ is described in~\cite{Mu}. 
The set of irreducible modules is described in Prop.~5.2 there, a result first 
proven in~\cite{Mi-rep}. In special cases, the fusion algebra and the
$S$-matrix and the intertwining algebra were determined in~\cite{La1,La2}. 

Let now $W$ be a self-dual FVOA. From Theorem~\ref{thVOAmilgram},
one gets $16\mid r$, Th.~2.19 of~\cite{DGH-virs} gives 
$\mathcal{D}=\mathcal{C}^{\perp}$, and one has ${\bf 1}^r\in \mathcal{D}$ for
the all-one-vector. 

\begin{theorem}[Relative Mass formula for FVOAs, \cite{Ho-virmass}]\label{relmass}
For $16\mid r$, one has
$$\sum_V   \frac{1}{|{\rm Aut}_{L_{1/2}^{\otimes r}(0)}(V)|}=
  \frac{1}{2^r\cdot r!}\cdot \sum_{k} 2^{\frac{k(k-1)}{2}+1}\,\sigma_k(r),$$
where the sum on the left hand side of the equation
runs over all  isomorphism classes (in the sense of extensions) $V$
of self-dual FVOAs of central charge $r/2$ and
$\sigma_k(r)$ is the number of codes $\mathcal{D}$ of length~$r$
with $\dim \mathcal{D}=k$, ${\bf 1}^r\in  \mathcal{D}$
and \hbox{${\rm wt}(\mathcal{D})\subset 8{\bf Z}$.}
\end{theorem}

{\bf Sketch of proof:} Since $V_\mathcal{C}$ is a simple current extension,
its modular  category $\mathcal{T}(V_\mathcal{C})$ can be expressed
completely in terms of the known modular category $\mathcal{T}(L_{1/2}^{\otimes r}(0))
=\times_{i=1}^r \mathcal{T}(L_{1/2}^{\otimes r}(0))$ and the code $\mathcal{C}$. 
For $V$ self-dual, all quantum dimensions of $\mathcal{T}(V_\mathcal{C})$ are~$1$,
so it is an abelian modular category $\mathcal{T}(A,q)$ with 
an explicitly given abelian group $A$ and quadratic form $q$.
One has a short exact sequence
$$\begin{array}{ccc}
{\bf F}_2^n/\mathcal{C}& \longrightarrow& A\\
 &  & \downarrow \\
 & & \mathcal{D}
\end{array}$$
and the self-dual extensions of $V_\mathcal{C}$ are in one-to-one correspondence
with sections $s: \mathcal{D}\longrightarrow A$ for which $s(\mathcal{D})\subset (A,q)$
is isotropic. Counting them gives the result.
\phantom{xxx}\hfill \qed

\smallskip 

Let us check the theorem for the first case $r=16$. In~\cite{GH-stabil}, 
the five FVOAs of central charge~$8$ together with their automorphism group have been determined.
For every $k=1$, $2$, $3$, $4$, $5$, there is an up to equivalence unique code $\mathcal{D}_k$
of dimension $k$ and a unique FVOA $V_k$ with code $\mathcal{D}(V_k)=\mathcal{D}_k$ 
and automorphism group 
${\rm Aut}_{L_{1/2}^{\otimes 16}(0)}(V_k)=G_\mathcal{C}^k.{\rm Aut}(\mathcal{D}_k)$.
By~\cite{GH-stabil}, main Theorem~2, one has
$|G_\mathcal{C}^k|=2^{15-k(k+1)/2+k}$ and
for the sum $\sum_V \frac{1}{|{\rm Aut}_{L_{1/2}^{\otimes 16}(0)}(V)|}$
we obtain:
$$\sum_{k=1}^5\frac{1}{|{\rm Aut}_{L_{1/2}^{\otimes 16}(0)}(V_k)|} =
 \frac{1}{2^{16}\cdot 16!}\, \sum_{k=1}^5 2^{1+k(k+1)/2-k}\cdot
\frac{16!}{|{\rm Aut}(\mathcal{D}_k)|},$$
in agreement with Theorem~\ref{relmass}.

The methods of the proof can be used to give a construction of
the Moonshine module $V^{\natural}$ as the self-dual FVOA with 
$\mathcal{C}=\mathcal{D}^{\perp}$ the lexicographic code of length~$48$ and minimal 
weight~$4$; see~\cite{DGH-virs}, Sect.~5 for the precise structure of the
Virasoro module decomposition.
The only input from VOA-theory which one needs is the construction of the Virasoro VOA 
of central charge $1/2$ and the structure of its intertwining algebra. 
All previous constructions use at some place the lattice vertex operator
algebra construction; cf.~\cite{FLM,Hua,Mia-Moonshine}.


\small

\providecommand{\bysame}{\leavevmode\hbox to3em{\hrulefill}\thinspace}

\end{document}